An Exploration Into the Collatz Conjecture with Changed Parameters

Raina Shrimali, Shouvik Ahmed Antu, and Miranda Jones

November 18, 2022

Willamette University

IDS-101-32-22/FA

*Abstract:*

Exploring the Collatz Conjecture and changing the expression from 3n + 1 to 5n + 1, we found patterns in different sets of numbers. Some numbers reduce to one (as stated in the Collatz Conjecture), some might escape to infinity, and some get stuck in repeating cycles. To further explore the patterns involved in Collatz-like expressions, we changed the expression to 3n+5 and found connections between 3n+1 and 3n+5.

*Introduction:*

The Collatz Conjecture is an interesting, arbitrary, and seemingly simple problem that mathematician Luther Collatz introduced to the world two years after receiving his doctorate in 1937 (O'Connor, Robertson). It has continued to stump the mathematical world since, as no one has been able to prove that it is either false or true. The conjecture goes as follows:

Take any natural number. If the integer is even, it is divided by two until it becomes odd. Once it is odd, it is plugged into the equation "3n + 1". According to Collatz, if these two operations are performed on each term derived from the starting integer, the end term will eventually reach one. The conjecture can be tested using the following piecewise function:

$$f(n) = \begin{cases} 3n+1 & : n \text{ is odd} \\ \frac{n}{2} & : n \text{ is even} \end{cases}, \text{ where n is a positive integer.}$$

However, the original Collatz Conjecture has been worked on numerous times, and the remaining work on actually proving the conjecture requires more resources and experience than our group has. Mathematicians have proven thus far that any number below $2^{68}$ plugged into the 3n+1 expression will always return to one (Cooke). Terence Tao made famous progress on the problem when he proved in 2019 that *most* numbers - 99% or more - will return to one, although he failed to prove the conjecture true in all cases (Hartnett). As we do not have a super computer, we decided to turn our research towards other interesting patterns, and landed on the modified expression "5n+1" instead, which has some unique and interesting properties. Our reason for choosing to multiply *n* by five instead of three is that this simple variation creates outcomes that deviate from the Collatz Conjecture. The odd numbers that are plugged into this equation don't all return to one; some get stuck in cycles and some repeat infinitely. Because of this, our primary research focus is the properties of the number loops as well as which numbers get stuck in loops and why.

*Mathematical Process:*

Our primary work on the modified version of the Collatz Conjecture involved testing numbers 1-100. The reason for testing this limited pool of numbers is that multiplying integers by five causes them to increase very quickly. The bigger the number we start with, the more chance that it will (seemingly) escape to infinity. Thus, we only tested 100 integers unless we found a pattern that needed further exploring. For instance, we found a pattern within multiples of 2 and 5 below 100, so we tested numbers greater than 100 to see if the pattern continued to show up.

Our research led us to see patterns almost immediately, with the help of our Python computer program. The interested reader may find the code in **Appendix A**. The program showed us whether or not these numbers repeat (infinitely or if they get stuck in a repeating cycle) without the need for doing the equations by hand. We made a Google Spreadsheet in which we recorded whether or not each number repeats, the number of digits in each number cycle, and whether or not the number eventually reaches one. We recorded this information in order to look for patterns. There was no specific pattern we had in mind; we looked for anything out of the ordinary that could eventually turn into a provable conjecture.

The process of searching for patterns went much quicker than we initially thought it would. While writing out numbers on a whiteboard, a pattern soon emerged within numbers $2^r * 5$ where ($r \in \mathbb{N}$). These numbers all have repeating cycles that do not return to one.

| $2^r * 5$ where $(r \in \mathbb{N})$ | Total numbers in repeating cycle: | Number cycle: |
|---|---|---|
| 5 (base number) | 10 | {26, 13, 66, 33, 166, 83, 416, 208, 104, 52} |
| 10 | 10 | {26, 13, 66, 33, 166, 83, 416, 208, 104, 52} |
| 20 | 10 | {26, 13, 66, 33, 166, 83, 416, 208, 104, 52} |
| 40 | 10 | {26, 13, 66, 33, 166, 83, 416, 208, 104, 52} |
| 80 | 10 | {26, 13, 66, 33, 166, 83, 416, 208, 104, 52} |

**Figure A**: shows the repeated cycle of numbers $2^r * 5$, where $(r \in \mathbb{N})$.

We found that some multiples of two and five have the same repeating number pattern that contains 10 different numbers. This is because each of these numbers turn into five after being divided by two. The amount of times each number is divided by two has no effect on the repeating number cycle. Five is also included in this group, which is why we had to name this group of numbers $2^r * 5$ where $(r \in \mathbb{N})$. Our theorem is:

For the piecewise function ,

$$f(n) = \begin{cases} n/2, & \text{when } n \text{ is even} \\ 5n + 1, & \text{when } n \text{ is odd} \end{cases}$$

when n is equal to $10 * 2^r$, the repeating cycle of the infinite loop starts after $(1 + r)$ calculations.

After making this theorem, we wanted to see what other number cycles exist. We. found that there are only two variations of numbers within the repeating cycles of numbers 1-100. This finding may later turn into a conjecture involving all positive integers once we test a higher quantity of numbers. Some non-multiples of two and five have cycles that include the same numbers listed in the table above and some have a completely different set of 10 numbers:

| Number: | Total numbers in repeating cycle: | Number cycle: |
|---|---|---|
| 13 | 10 | {66, 33, 166, 83, 416, 208, 104, 52, 26, 13} |
| 17 | 10 | {86, 43, 216, 108, 54, 27, 136, 68, 34, 17} |
| 26 | 10 | {13, 66, 33, 166, 83, 416, 208, 104, 52, 26} |
| 27 | 10 | {136, 68, 34, 17, 86, 43, 216, 108, 54, 27} |
| 33 | 10 | {166, 83, 416, 208, 104, 52, 26, 13, 66, 33} |
| 34 | 10 | {17, 86, 43, 216, 108, 54, 27, 136, 68, 34} |
| 43 | 10 | {216, 108, 54, 27, 136, 68, 34, 17, 86, 43} |

**Figure B**: shows the number cycles of the first six numbers with repeating cycles. The numbers within the cycles are written in the order they appear when plugged into 5n + 1.

▨ = first cycle variation (seen above in multiples of both two and five)

▨ = second cycle variation

Numbers in order from smallest to biggest:

13, 26, 33, 52, 66, 83, 104, 166, 208, 416

17, 27, 34, 43, 54, 68, 86, 108, 136, 216

Numbers in the sequence they appear (starting at 13 and 17):

13, 66, 33, 166, 83, 416, 208, 104, 52, 26

17, 86, 43, 216, 108, 54, 27, 136, 68, 34

Cycles combined in order from smallest to biggest (to show the diversity of numbers within the cycles):

13, 17, 26, 27, 33, 34, 43, 52, 54, 66, 68, 83, 86, 104, 108, 136, 166, 208, 216, 416

In comparison, the 5 vs 13 cycle in sequence:

5: 26, 13, 66, 33, 166, 83, 416, 208, 104, 52

13: 66, 33, 166, 83, 416, 208, 104, 52, 26

▇ = odd number*

* The repeating number cycles are very puzzling. Every odd number ends with either a three or a seven. Our group has not yet discovered the reason behind this strange pattern.

    Our group was getting a bit hung up on deciphering a pattern within the number cycles. Why do these two cycles continue to appear? Why do the numbers contained within them repeat infinitely? Are there more cycles that we haven't found yet because we haven't been able to test a wider range of numbers? As we continue to collect data and test numbers, the answers to these questions will gradually appear. In the meantime, we created a simple but effective theorem for avoiding confusion in the face of repeating number cycles:

    If you take a number from either cycle, *m,* and multiply it by $2^r (r \in \mathbb{N})$, the resulting number will have the same repeating number cycle as *m* when applying the piecewise function

$$f(n) = \begin{cases} n/2, & \text{when } n \text{ is even} \\ 5n + 1, & \text{when } n \text{ is odd} \end{cases}$$

This is due to the fact that the number (m * $2^r$) must be divided by two before being entered into the equation 5n + 1, and these divisions will result in the number eventually returning to *m*. The number of times the number is divided by two does not affect the repeating cycle in any way. This theorem has limitations because we haven't gotten around to testing numbers containing more than four digits. However, for our current purposes and available resources, this theorem works while we stick to testing out the smaller numbers.

After exploring 5n + 1, we changed the expression to 3n+5 because it contains the primes 3 and 5 and also because the numbers that are generated are close to the numbers that are generated when 3n+1 is used. However, unlike 5n + 1 which may cause numbers to escape to infinity, 3n + 5 seems to only produce repeating loops. There are four repeating cycles that we have found so far.

| Number: | Total numbers in repeating cycle: | Number cycle: |
|---|---|---|
| 1 | 4 | 8, 4, 2, 1 |
| 3 | 8 | 38, 19, 62, 31, 98, 49, 152, 76 |
| 5 | 3 | 20, 10, 5 |
| 23 | 8 | 37, 116, 58, 29, 92, 46, 23, 74 |

**Figure C**: shows the four different number cycles that have shown up in our computations so far. The numbers in the left column are chosen because they are the first numbers that create each different number cycle.

After making a new code for 3n + 5, my group followed the same process that we used to explore the expression 5n + 1 and recorded the number cycles of the first 100 natural numbers in a Google spreadsheet. At first it seemed like the repeating cycles occurred at random. However, we developed a new conjecture after noticing a pattern:

Multiples of 5 will create the repeating cycle *20, 10, 5* when applying the piecewise function

$$f(n) = \begin{cases} n/2, \text{ when } n \text{ is even} \\ 3n + 5, \text{ when } n \text{ is odd} \end{cases}$$

| Number: | Total numbers in repeating cycle: | Number cycle: |
|---|---|---|
| 5 | 3 | 20, 10, 5 |
| 15 | 3 | 20, 10, 5 |
| 25 | 3 | 20, 10, 5 |
| 35 | 3 | 20, 10, 5 |

**Figure D**: shows the first four multiples of five and their repeating cycles.

| Number: | Total numbers in repeating cycle: | Number cycle: |
|---|---|---|
| 225 | 3 | 20, 10, 5 |
| 17585 | 3 | 20, 10, 5 |

| 3698450 | 3 | 20, 10, 5 |

**Figure E**: shows three random numbers ending in five or zero, and their repeating cycles. This is done to show that not just two digit multiples of five have the same repeating cycle.

Now notice, the numbers in 20, 10 and 5 are exactly 5 times of the numbers 4, 2 and 1. Going back to the Collatz Conjecture, all the numbers reduce following the pattern 4, 2, and 1, and hence this is a direct connection between the Collatz Conjecture and our new expression. Notice the following diagram, the relation between the two expressions:

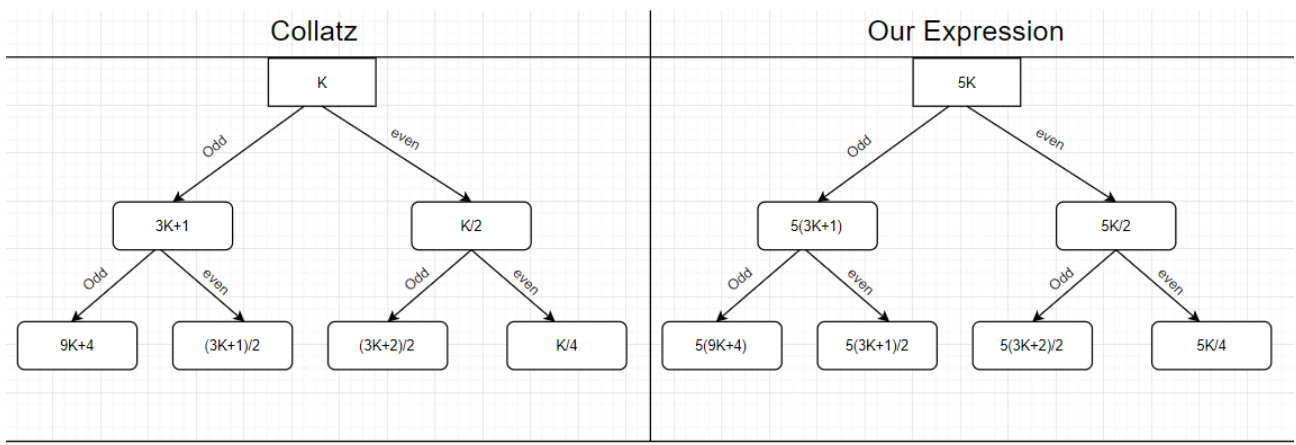

We are considering any positive integer and the two possibilities a number can be, that is, even or odd. Taking any arbitrary number, $k$, and plugging it into the piecewise function that represents the Collatz Conjecture, then plugging 5 times of that number, $5k$, into our new expression, we can see that the numbers are 5 times of the output that we get from the Collatz piecewise function.

Since we are considering both cases of being even and odd, we have proved that if Collatz conjecture holds:

> Given any positive arbitrary integer K for the functions X, Y below, the computations will result as *y(k)=5(x(k))*

$$X(k) = \begin{cases} k/2, & \text{when } k \text{ is even} \\ 3k+1, & \text{when } k \text{ is odd} \end{cases}, \quad Y(5k) = \begin{cases} k/2, & \text{when } k \text{ is even} \\ 3k+5, & \text{when } k \text{ is odd} \end{cases}$$

Hence, the outputs from our expression 3n + 5 correspond and have a connected pattern with Collatz's original expression.

*Conclusion:*

The Collatz Conjecture is a well known conjecture that has been thoroughly studied for decades by many exceptional mathematicians. Our slight variations on this conjecture have led to interesting discoveries. Working with 5n + 1, we found that some numbers act like how Collatz predicted, some go into repeating loops and some might actually escape to infinity. Working with 3n + 5, we noticed that none of the numbers seem like they escape to infinity, all of them have some repeating cycle, and, more interestingly, only 4 different repeating patterns show up when plugging numbers into the piecewise function. The connection between our 3n + 5 expression and Collatz's original expression is a unique discovery that opens the door for further exploration.

**Appendix A:**

Python Code:

```
x= [Any Natural Number]
listofnum=[ ]
while x !=1  and len(l) != 100000:
   if x%2==0:
      x=x/2
      l.append(int(x)
   else:
      x=(5*x+1)
      l.append(int(x))
print(listofnum)
```